\documentclass [a4paper,14pt ]{extarticle}
\usepackage[cp1251]{inputenc}
\usepackage[russian]{babel}
\usepackage[dvips,pdftex]{graphicx,graphics}
\usepackage{amsfonts}
\usepackage[dvipsnames,pdftex]{color}

\begin{document}
 \baselineskip=16pt
\clearpage

\newcommand{\N}{\mathbb{N}}
\newcommand{\qq}{\qquad}
\newcommand{\q}{\quad}

\newcommand{\si}{\sigma}
\newcommand{\la}{\lambda}
\newcommand{\al}{\alpha}
\newcommand{\be}{{\beta}}

\vspace{4ex}

\baselineskip=16pt

 \setcounter{page}{1}

{\ }

\vspace{3ex} 

 MSC-2010-class: 11A25 (Primary); 11B83 (Secondary)

\begin{center}

\vspace{3ex}

{\ } {\bf  \Large One-step  $G$-unimprovable numbers\footnote{\ \ \ This work
 was supported
by the grant of Russian Foundation of Fundamental Research\\
$\qq {\ } {\qq } {\ } $ (project \# $14-01-00684.$)}
}

\vspace{2ex}

{\bf \large by Gennadiy Kalyabin\footnote{\q Samara State Technical University, Russia; \ gennadiy.kalyabin@gmail.com}}

\end{center}
\vspace{1ex} 

{\it Abstract:} 
The infinitude is established  of the set ${\bf U_1}$ 
of all integers $N>5$ whose Gronwall numbers $G(N)$ do 
not increase when replacing $N$ by $N/q$ or $Np$, where $q, p$ are primes, $q\ | N$.

\vspace{1ex}


\vspace{1ex}

Bibliography: 6 items

\vspace{3ex}

{\bf 1. \ Notations and problem setting}

\vspace{2ex}

 As usually let $\log x$ and  $\sigma(n)$ stand (resp.) for the natural logarithm of a positive $x$ and the sum of all divisors of a positive integer $n$.  In 1913  T. Gronwall established  [1] the limit relationship which involves the Euler-Masceroni constant $\gamma:=\lim_{n\to\infty}\left(\sum_{k=1}^n \frac{1}{k} -\log n\right)=0.577215..$: 
$$\limsup_{n\to\infty} \frac{\sigma(n)}{n \log \log n}=e^{\gamma}=1.78107..
\eqno{(1.1)}
$$

S. Ramanujan has noticed (in 1915, the first publication in 1997 [2]), that: \
{\it if the Riemann Hypothesis 
on non-trivial zeros of $\zeta(s)$ holds true, then in addition to (1) for all sufficiently large  $n$ the {\bf strict} inequality
$$ {\sigma(n)}<e^{\gamma}{n \log \log n}; \q n>n_0,
\eqno{(1.2)}
$$
is fulfilled.}   70 years later G. Robin [3] proved a  paramount  assertion  which in a sense is inverse to Ramanujan's result, namely:\  

\vspace{1ex}
{\it if  (1.2) holds true for all integers  $n>5040$, then RH is valid. }

\vspace{2ex}
Detailed discussion of historical aspects and adjacent questions may be found in a remarkable  Caveney-Nicolas-Sondow paper [4], (which in fact was the initial point for this author's research) where it is proved that RRI is
in turn equivalent to the statement:

\vspace*{\fill}

\clearpage

\vspace{1ex}

 {\it Every integer $N>5$ is $G$-improvable, i. e.  either:} 
\vspace{1ex}

$\mathbb{\bf I}^/:$ {\it there is a prime } $q\ | \ N, \hbox{\it \ such that  }G(N/q) > G(N),$ \q {\it or:}

\vspace{1ex}

$  \mathbb{\bf I}^{\times}:$ {\it there is an integer} $ a>1$ \ {\it such that} \ $G(Na)> G(N)$. 

\vspace{1ex}

Note that the numbers 3, 4 and 5  {\bf are not}  {\it $G$-improvable.}
 
\vspace{1ex}

We will consider a new class ${\bf U}_1$  of all such integers $N>5$ which cannot be $G$-improved neither by multiplication 
nor by division by any {\bf single} prime.

\vspace{1ex}

{\ D\ e\ f\ i\ n\ i\ t\ i\ o\ n}. \ {\it An integer $N>5$ is called 1-step $G$-umimprovable ($N\in {\bf U_1}$) if and only if the following} 
{\it two conditions hold:} 

\vspace{1ex}

$ {\bf U_1^/}:$ {\it for any prime } $q\ | \ N, \hbox{\it \ one has  }G(N/q) \le  G({N}),$ \q 
{\it and:}


\vspace{1ex}

$  {\bf U_1^{\times}}:$ {\it for any prime  $ p$ \  one has} \ $G(Np)\le G(N)$. 

\vspace{1ex}

 R\ e\ m\ a\ r\ k\q 1. \ The condition 
 ${\bf U_1^/}\ $ (in [4], S. 5 it was studied under
 the name GA1) is exactly the negation
 of ${\bf I}^/$ whereas  ${\bf U_1^{\times}}$ is essentially weaker than the negation of
 ${\bf I^{\times}}$. Thus the CNSH is equivalent 
to ${\bf U_1}\subset {\bf I^{\times}}$. 

The purpose of this paper is to establish the infinitude of   ${\bf U_1}$ and to construct the explicit algorithm
which successively calculates all elements of this class the least of them being equal 
$$N^*_1= 2^5\cdot 3^3 \cdot 5^2 
\cdot 7 \cdot 11 \cdot 13  \cdot 17 \cdot 19 \cdot 23 
= 160 \ 626\ 866\ 400;\ G(N_1^*)=1.7374...
\eqno{(1.3)}
$$

The author believes this approach to be helpful for the proof of the CNSH.

\vspace{2ex}

{\bf Preliminary results}
\vspace{1ex}

Let  $\N$ be the set of all positive integers, $\N_0:=\N\cup \{0\}, \ (a,b):=\gcd(a,b) \ -$ the greatest common divisor of  $a,b \in \N$;
 $\mathbb{P}:= \{p_k\}_1^{\infty}  $ -- increasing sequence of all primes; the notation  
$$N\ \| \ p^{\al},\ p\in \mathbb{P},\ \al \in \N_0,  \hbox{ means that  } N=p^{\al}m, \ (m,p)=1.
\eqno{(2.1)}
$$
 \vspace{1ex}

First we will study how the Gronwall number changes if one replaces  $N$ by $Np,\  p\in \mathbb{P}$. Let us denote 
$$ \la \equiv \la(p,\al):= \frac{p^{\al+2} -1}{p^{\al+2} -p}=1+\frac{1}{p+p^2 + \dots +p^{\al+1}};
\eqno{(2.2)} 
$$
and let  $\xi:=\xi(p, \alpha) $ be the unique positive root of the equation 

\vspace{2ex}

$$ \xi^{\lambda(p, \alpha)}= \xi+\log p 
\iff \xi^{\lambda(p, \alpha)-1}
= 1+\frac{\log p}{\xi}
$$
$$
\iff \ \log \xi =(p+p^2+\dots +p^{\al+1})
\log\left(1 + {\log p \over \xi}\right).
\eqno{(2.3)}
$$

 \vspace{1ex}

L\ e\ m\ m\ a\ \  1. \ {\it Let $N>1, \ N \ \| \ p^{\al} 
$;  then three following conditions\\ are equivalent  (comp. the definition of ${\bf U_1^{\times}}$): }
$$\hbox{\bf (i)}\ G(Np) >  G(N),\ \hbox{\bf (ii)} (\log N)^{\la} > \log N+\log p, 
\ \hbox{\bf (iii)} \log N > \xi(p, \alpha).
\eqno{(2.4)} 
$$

$\triangle$ By virtue of the classical  number theory  formula (cf.  [5], Ch 1, \S\ 2d):
  $$\ (p, m)=1\ \Rightarrow \  \si(p^{\al}m)=\frac{p^{\al+1}-1}{p-1}\ \si(m);
\eqno{(2.5)}
$$

whence the identities follow:
$$\frac{G(p^{\al+1}m)}{G(p^{\al}m)}=
\frac{p^{\al+2}-1}{p-1}\cdot \frac{p-1}{p^{\al+1}-1}\cdot \frac{p^{\al}  \log \log (p^{\al} m )}{p^{\al+1} \log \log (p^{\al+1} m )}
$$
$$=  \frac{\la   \log   ( \log N)}{\log( \log N+ \log p )}
=  \frac{ \log (\log N)^{\la }}
{\log( \log N+ \log p )}.
\eqno{(2.6)}
$$

Therefore the inequality  $G(Np) > G(N)$ is
equivalent to
$$ \frac{G(p^{\al+1}m)}{G(p^{\al}m)}  > 1 \iff (\log p^{\al}m)^{\la} 
> \log p^{\al}m + \log p;
\eqno{(2.7)}
$$
this proves the equivalence of (i) and (ii) in (2.4). It remains to note that since the function
 $\la(p, \al)$ {\it decreases } as either $p$ or 
 $\al$ increase (cf. (2.2)), the root  $\xi(p, \al)$ of (2.3) is an {\it increasing } function of both  $p$ and
 $\al \q \Box$


\vspace{1ex}

L e m m a \  2. {\it Let the exponent ${\al}$ in Lemma 1 be positive; \ then the following three conditions are equivalent (cp. $\ {\bf U^/}$)}:
$$  \hbox{\bf (i)  }  \ G(N/p) \le G(N), \q \hbox{\bf (ii)  }  \ (\log N - \log p)^{\la(p, \al -1)} \ge \log N,
$$
$$ \hbox{\bf (iii)  }  \ 
 \log N \ge \log p+\xi(p, \al -1).
\eqno{(2.8)} 
$$


$\triangle$ \ It suffices to apply Lemma 1 to the numbers 
$p$ and $\tilde{N}:=N/p \ \|\ p^{\al -1} \ \Box$

\vspace{1ex}

R e m a r k \ 2. The assertion of Lemma 2 
in somewhat different form was given in 
[4, Proposition 15]. 
It was also observed that the integer
$$ \nu:=2^4\cdot 3^3 \cdot 5^2 \cdot 7 \cdot 11
\cdot 13 \cdot 17=183\ 783\ 600
\eqno{(2.9)} 
$$
is the least $N>4$ such that for all  primes
$q \ | \ N,$ the inequality $G(N/q)<G(N)$ is fufilled. However 
$G(19\nu)=1.7238.. >G(\nu)=1.7175..$ and thus $\nu \not\in
 {\bf U_1}$.

\vspace{1ex}

Joining the two Lemmas yields the assertion helpful 
in  the sequel. 

\vspace{1ex}

P r o p o s i t i o n \ 1.\ 
{\it \ Let $N \ \| \ p^{\alpha}, \ \alpha>0$;  then following two}  

{\it  properties are equivalent}:

\vspace*{\fill}

\clearpage

 $${\bf (i)  }\ N \ { is }\ G_p\hbox{-}unimprovable, 
\hbox{ i. e. } 
G(N)\le \min (G(N/p), G(Np)) , 
$$
$$ {\bf (ii)  } \log N\in \Delta_{p, \alpha}
:=       [ \xi(p, \alpha -1) +  \log p,\ \xi(p, \alpha)] .
\eqno{(2.10)}
$$

{\it For fixed prime $p$ the segments $\Delta_{p, \alpha}, \ \al\in\N$ are nonempty and disjoint; in the case when $\log N$ belongs to the junction interval  $I_{\alpha}:=(\xi(p, \alpha), 
\xi(p, \alpha) +\log p)$, the number  $Np$  is
 $ G_p$-unimprovable, i. e. } $G(Np)\le \min(G(N), G(Np^2))$.

\vspace{1ex}

We also  need the bilateral estimates for the 
roots of equations (2.3) (more sharp as $\al=0$)
which are obtained by standard Analysis techniques.

\vspace{1ex}

L e m m a\ \ 3.\ {\it For all $\alpha\in {\mathbb{N}}_0, \ p\in \mathbb{P}$ the following inequalities are fulfilled: }
$$ \hbox{\bf (i) }\   p \ -  \ { \log p}  <\xi(p,0)< \ p; $$
$$
\hbox{\bf (ii) } \   \frac{p^{\alpha+1}}{ \alpha +1}
 \ <\xi(p,\alpha)\ <  3\ \frac{ p^{\alpha+1}}
{ \alpha +1}, \quad \alpha\in \mathbb{N}. 
\eqno{(2.11)}
$$

$\triangle$ \ {\bf (i) } The determining equation 
(2.3) $\xi^{1+{1\over p}}
=\xi + \log p$ for  $\xi(p,0)$ may be rewritten in the form:
 $$ f(\xi):={\log \xi \over p} -
\log\left(1+{\log p \over \xi}\right)=0,
\eqno{(2.12)} 
$$
where $f(\xi)$ is a strictly increasing 
function   for $\xi>0$.  
From the first inequality
$$ \hbox{ (i) }  \log(1+z)<z, \ (z>-1, \ z\not=0), \q 
$$
$$\hbox{ (ii) }  \log (1+z)>z-0.5 z^2 \q (z>0) 
\eqno{(2.13)}
$$ 
and (2.12) it follows that $f(p)>0$.  Further 
 for  $\xi_1:=p-\log p$ one obtains analogously 
$$  f(\xi_1) = {\log (p-\log p)\over p}
-\log\left(1+{\log p\over p-\log p}\right)
$$
$$  = - \log p 
+\left(1+{1\over p}\right) \log (p-\log p)
 ={\log p\over p}+\left(1+{1\over p}\right) \log \left(1-{\log p\over p}\right)
$$
$$ < {\log p \over p}
-\left(1+{1\over p}\right){\log p\over p}
= -{\log p\over p^2}< 0 \qq \Box
\eqno{(2.14)}
$$

 \vspace{2ex}

$\triangle$ \ {\bf (ii)} Now let $\al>0$; by  introducing two new notations: 
 $$w=w(p, \al):=1+{1\over p}+ {1\over p^2}+\dots  +{1\over p^{\al}}; \q t:={(\al+1)\xi \over p^{\al+1}},\ \hbox{ i. e.  }\ \xi={\, t\, p^{\al+1}\over \al +1}; 
\eqno{(2.15)} 
$$
one can rewrite the determining equation (2.3)
in the following equivalent form:
$$ g(t):=\log t -\log (\al+1) + (\al+1) \log p
$$
$$ - p^{\al +1}\ w\ \log\left(1+{(\al+1) \log p \over t\, p^{\al+1}}\right)=0. 
\eqno{(2.16)} 
$$

For  $p, \al$ being fixed the function  $g(t)$  
is strictly increasing; therefore the assertion 
 {\bf(ii)} of Lemma 3 is equivalent to  two inequalities:  $g(3)>0,\ g(1)<0$.

 \vspace{1ex}

To prove the first of them we again  apply 
 (2.12){\bf (i)}, and taking into account that
$w<2$, we obtain
$$ g(3)> \log 3 -\log (\al +1)  +  (\al +1) \log p
 -  {w\over 3} (\al +1) \log p  
 \eqno{(2.17)}
 $$
$$ > \log 3 - \log (\al +1) +\frac{\al +1}{3}\log p
\ge  1+\log \log p \ge 1 + \log \log 2 = 0.633..
 $$

 \vspace{1ex}

To  demostrate that $ g(1)< 0$, we take use of
 the inequality  (2.12){\bf (ii)}:
$$ g(1)<  -\log (\al +1) +  (\al +1) \log p
 -  w (\al +1) \log p\  \left( 1  
- \frac{(\al +1) \log p}{ 2\, p^{\al +1}}\right)
$$
$$ = - \log (\al +1)  +  
\left( 1  - w\left(1 -\frac{(\al +1)  \log p}{ 2\, 
p^{\al +1 }}\right)\right)  (\al +1)  \log p
\eqno{(2.18)} 
$$
$$ < - \log (\al +1)   +  
\left( 1  - \left(1+{1\over p} \right) 
\left(1 -\frac{(\al +1)  \log p}{  2\, p^{\al +1 }}\right)\right)  (\al +1)  \log p <0.
$$

Here we have taken into account the relationships valid for all $p\ge 2,\ \al \ge 1$  
 $$w\ge 1+{1\over p}; \q \frac{(\al +1) \log p}
{ 2 p^{\al +1 }} \le {\log p\over p^2}\q \Longrightarrow 
\ \left(1+{1\over p}\right)\left(1-{\log p\over p^2}\right)>1,
\eqno{(2.19)}
$$
whence it follows that the large bracket in the last line (2.18) is also negative 
$\  \Box$

\vspace{1ex}

\vspace{1ex}

The rather simple assertion concerning non-negative sequences will be also frequently helpful for us.

\vspace{1ex}

L e m m a\ \ 4.\ {\it  
Let  $A>0, B>0;$ denote by
$$ \tau^*=\tau^*(A, B):= A+0.5B(B+\sqrt{B^2+4A}) \ 
\eqno{(2.20)}
$$ 
the  root of the equation  $\tau=A+B\sqrt{\tau}$;\ suppose that $ a_n\ge 0,$
$ a_{n+1}\le A+B\sqrt{a_n}$ 
  $\forall n\in {\N}_0$;
\ then if $a_0 \le \tau^*$,  
then  $a_n\le \tau^*$  for all} $ n$.

$\triangle$ It suffices to consider the auxillary sequence $b_n\ge a_n$ defined recursively as 
$b_0:=a_0, \ b_{n+1}= A+B\sqrt{b_n}$ 
which increases to its limit $\tau^*\ \Box$

\vspace*{\fill}
\clearpage

{\bf 3. \  Structure and properties of numbers from the class  ${\bf U_1}$ }\ 

\vspace{1ex}

We will use the first and the second Chebyshev functions (cf. [6],  \ 3.2, p. 104)
 $$\theta(x):= \sum_{p\le x} \log p; \quad
\psi(x):= \sum_{p^k\le x} \log p = \theta(x) +\theta(x^{1/2})+ \theta(x^{1/3})+\dots,
\eqno{(3.1)}
$$
as well as Chebyshev products
 $T(p_n):= \exp(\theta(p_n)) = p_1\cdot p_2\cdot \dots\cdot p_n$. 

Let   $P(N)  $  stand for the greatest  prime factor of an integer $N>1$.

\vspace{1ex}
  
T h e o r e m \  1.
 {\it \ Let the integer $ N>5 $ have the canonical factorization }
 $$ N =p_1^{\alpha_1}\ p_2^{\alpha_2}\ p_3^{\alpha_3} \dots p_k^{\alpha_k}, \quad \alpha_k>0, \ \hbox{ i. e. } P(N)=p_k. 
\eqno{(3.2)}
$$

 {\it Then} \ 
{\bf (A)}:  {\it the belonging  $  N\in {\bf U_1}$ is equivalent to the set of inequalities:}
 $$ {\bf (i)}\  \xi(p_j,\alpha_j -1) 
+\log p_j  \le \log N \q \forall j\le k \hbox{ such that}
\ \al_j>0; 
$$
$$ {\bf (ii)}\q \log N \le \xi(p_i,\alpha_i) \q 
\forall i \le k; \q {\bf (iii)}\ \log N \le \xi(p_{k+1}, 0).
\eqno{(3.3)}
$$

\vspace{1ex}

${\bf (B_1)}$:\ {\it the exponents $\alpha_j$ do not increase: \ $\alpha_1\ge \alpha_2\ge     \dots \ge \alpha_k$ and}  $\alpha_k=1$.

\vspace{1ex}

$\triangle$ \ 
 The equivalence of conditions 
$ N\in {\bf U_1}$ and 
  ({3.3}){\bf (i)-(iii)} follows immediately from
definitions of  ${\bf U_1^/}$ and  ${\bf U_1^{\times}}$ и and assertions of lemmas 2 and 1 with taking into account that $\xi(p_n, 0) > \xi(p_{k+1}, 0)$ for all  $n>k+1$. 

Further assuming that  $\alpha_i<\alpha_j$ for some   $i, j,\ 1\le i<j \le k$, one obtains bearing in mind the
monotone increase of  $\xi(p,\alpha)$ with respect
to both  $p$ and  $\alpha$ 
$$ \xi(p_j, \alpha_j-1)+\log p_j \ge \xi(p_j, \alpha_i)
+\log p_j  >  \xi(p_i, \alpha_i) ,
$$ 
and thus for all  $N>5$  the violation the exponents monotonicity would lead to the incompatabity of
inequalities (3.3){\bf (i)} and  (3.3){\bf (ii)} for
$\log N$.

 At last, if one supposes 
$\alpha_k\ge 2$, then from  
 (3.3){\bf (i)}, (2.11){\bf (ii)}   it follows  that
$$ \log N\ge \xi(p_k, 1)+\log p_k>0.5 p_k^2+
\log p_k >p_{k+1}> \xi(p_{k+1}, 0),
\eqno{(3.4)}
$$
which is in contradiction with (3.3){\bf (iii)} for 
$k\ge 2$\q  $\Box$

\vspace{2ex}

{\bf 4.\ Locally $G$-extremal numbers}

\vspace{1ex}

We will start with the algorithm which generates the uniquely defined sequence of positive integers
 $V_k$, possessing {\it some part of properties,} prescribed by Theorem 1, and then we will prove
that the class ${\bf U_1}$ of all one-step $G$-unimprovable numbers, i. e. such that
satisfy {\it all } inequalities (3.3),  is its infinite subsequence 
$\{V_{k_m}\}=:\{N^*_m\}$, with indices  
$k_1<k_2 <\dots < k_m< \dots $  also constructively calculated.

\vspace*{\fill}
\clearpage

T h e o r e m \ 2.\q {\bf (I)}\ {\it For any fixed  $k\ge 4$  there exists the least (and thus unique)
integer  $V_k:=p_1^{\alpha_{1,k}}\dots p_k^{\alpha_{k,k}}$, where the exponents
 $\alpha_{j,k}, 1\le j \le k$ satisfy the condition ${\bf (B_1)}$ of Theorem 1, and 
the inequalities are fulfilled:} 
$$\xi(p_j,{\alpha}_{j, k} -1) +\log p_j \ 
 \le \log V_k \le  \xi(p_j, {\alpha}_{j, k} ) ; 
\qquad 1\le \  j\ < k;  
\eqno{(4.1)}
$$

\vspace{1ex}

{\bf (II)}\ {\it If some other integer 
$\tilde{V}_k:=p_1^{\tilde{\alpha}_{1,k}}\dots p_k^{\tilde{\alpha}_{k,k}}$ 
possesses the same properties }

{\it \qq \ as $V _k$ (with replacing  $V_k$ by   $\tilde{V}_k$ and
$\alpha_{j,k}$ by   $\tilde{\alpha}_{j,k}$), 
then}  $V_k \ | \ \tilde{V}_k$.

\vspace{1ex}

 {\bf (III)}\ {\it Moreover the relationship
$ \log V_k - \theta(p_k)=C_k\sqrt{p_k}$   holds}

\qq \q {\it   where} $\ 0.5<C_k<3$. 



\vspace{1ex}

$\triangle$ \ {\it Step  1.}\  According to Lemmas 1 and  2 the left and the right 
inequalities (4.1) are equivalent (respectively) to the fulfilment for the number
 ${V}_k$ of the conditions  ${\bf U_1^/} $ and  ${\bf U_1^{\times}}$ in which the primes
 $q, p$ {\it are subject to the restrictions}
 $q, p <p_k$; in other words in definition of  ${\bf U_1} $ the transitions
$V_k \to V_k/p_k$ and  $V_k \to V_kp_n, \ n\ge k$ are excluded.

\vspace{1ex}

Let us put  $ Y_{\ k}^{(0)} =T(p_k) :=p_1\cdot  \dots \cdot p_k$ (cf. (3.1)), 
and then  define inductively for  $s\in \mathbb{N}$:
$${\beta}_{j, k, s}:= \max \{ \beta\in \mathbb{N}:\xi(p_j,\beta -1)+\log p_j \le Y_{\ k}^{(s-1)}\};\ j<k;
$$
$$
{\beta}_{k, k, s}:= 1; \qquad  Y_{\ k}^{(s)} :=  \prod_{j=1}^k \ p_{\ j}^{\beta_{j, k, s}}.
\eqno{(4.2)}
$$

For  $k\ge 4 \ $ one has $\log {Y}_k^{(0)}=\theta(p_k)\ge 
\log 210=5.347.. > \xi(2, 1)+\log 2=2.56.. \ ; $  therefore  by virtue of Lemma 2 (cf. (2.8){\bf (iii)} the relationship ${\beta}_{1, k, 1}\ge 2$ holds, 
and hence   
${Y} _{\  k}^{(1)} >{Y}_{\ k}^{(0)} $.
The increase of function  $\xi(p, \alpha)$ with respect to both $p$ and $\alpha$ implies the
condition  ${\bf (B_1)}$ of Theorem 1 for the exponents  
$\{ {\beta}_{j, k, s}\}_{s=0}^{\infty}, \\ j<k; \ k, s$ being fixed. 

From the first line of defining formula (4.2) it follows that  ${\beta}_{j, k, s} \ge {\beta}_{j, k , s-1}$, 
and thus  ${Y}_{\ k}^{(s)}\ge {Y}_{\ k}^{(s-1)}$. Moreover, by virtue of the same formula
one obtains for all $j\in \{1, 2, \dots k-1\},\ s\in \mathbb{N}$:
$$\xi(p_j,{\beta}_{j, k, s}-1 ) +\log p_j 
\  \le \log {Y}_{\ k}^{(s-1)}  < \xi(p_j, {\beta}_{j, k, s})+\log p_j. 
\eqno{(4.3)}
$$

\vspace{1ex}

{\it Step 2.} \  Further, let us denote  $q_m^{(1)}:=\max\{p_j: \beta_{j, k, 1}\ge m\}, \ m=1, 2, \dots , \beta_{1, k, 1}$. It is clear  that $q_1^{(1)}=p_k> q_2^{(1)}\ge \dots \ge q_{\beta_{1, k, 1}}^{(1)}=2$. 
From the estimates  (2.11)${\bf (ii)}$ one may conclude that for all sufficiently large  $k$ 
(for small  $k$ the numbers  $C_k$ in  ${\bf (III)}$ are calculated directly) and for all  
 $m>1$ there holds $ q_m^{(1)}=\tilde{C}_{k,m}^{(1)}p_k^{1/m}, \ 0.7<\tilde{C}_{k, m}^{(1)}<2; \ \beta_{1, k, 1}^{(1)}< 1.5 \log \log {Y}_{\ k}^{(0)}$, and hence:

\vspace*{\fill}
\clearpage

$$ \log Y_{\ k}^{(1)}= \log(Y_{\ k}^{(0)})
+C_k^{(1)}\sqrt{Y_{\ k}^{(0)}}
 +O(p_k^{1/3} \log \log Y_{\ k}^{(0)})< 
\theta(p_k)+2\sqrt{Y_{\ k}^{(0)}}.  
\eqno{(4.4)}
$$

Going on  to argue in the same way as in Step 1 (with replacing the upper index $(1)$ by $(s)$), 
one obtains $ \log Y_{\ k}^{(s)} < \log Y_{\ k}^{(0)}+ 2 \sqrt{\log Y_{\ k}^{(s-1)}}$.

According to Lemma 4  these relationships imply the upper estimate: 
$\log Y_{\ k}^{(s)} < \theta(p_k)+2 + 2\sqrt{1+\theta(p_k)}$ for all  $s\in\mathbb{N}$, 
and because the sequence $\{ Y_{\ k}^{(s)}\}_{s=0}^{\infty}$ does not decrease and  is
{\it integer-valued}, it necessarily stabilizes beginning with certain number  ${s}_0:=s_0(k)$, 
and then,  by virtue of  definition (??), beginning with  the same number the exponents sequences 
  ${\beta}_{j,k,s}$  stabilize as well, and the inequalities 
$\xi(p_j,{\beta}_{j, k}-1 ) +\log p_j \  \le 
\log {Y}_{ k} < \xi(p_j, {\beta}_{j, k})+\log p_j $  hold true.

\vspace{1ex}

{\it Step  3.}\ Now introduce the set of indices 
$E_0:=\{j<k: \log Y_k>\xi(p_j, \beta_{j,k}\}$.
If it will turn out to be empty then all conditions of Theorem 2 are fulfilled
for  $V_k^{(0)}:=Y_k$. Otherwise, i. e. if $\log Y_k\in I_{p_j, \beta_{j,k}}$ 
for some  $j<k$ (cf. (2.10) ff), we put $\alpha_{j,k}^{(0)}:=\beta_{j,k},$ and
in accordance with Proposition 1 define recurrently for $s\in \mathbb{N}$:
 $$\alpha_{j,k}^{(s)}:=\beta_{j,k}+1,\ j\in E_{s-1}; 
\qquad  \alpha_{j,k}^{(s)}:=\beta_{j,k},
\ j\not\in E_{s-1};
$$
$$
 V_k^{(s)}:=\prod_{j=1}^k p_j^{\alpha^{(s)}_{j,k}}, 
\quad E_s:=\{j<k:  \log V_k^{(s)}>\xi(p_j, \beta_{j,k}) \}.
\eqno{(4.5)}
$$

It's clear that $V_k^{(s)}  \ge V_k^{(s-1)}, \ E_{s-1}\subset E_{s},\ \alpha_{j,k}^{(s)}\ge \alpha_{j,k}^{(s-1)}$. It should be emphasized
 that {\it no integer } $j\in E_{s-1}$
may {\it repeatedly} arise in  $E_{s} \setminus E_{s-1}$, because this would imply that 
$\sum\{\log p: p\in E_{s-1}\} > \xi(p_j, \beta_{j,k}+1)   - \xi(p_j, \beta_{j,k})-\log p_j$ 
what is impossible by virtue of estimates  
adduced earlier.

Now using the arguments analogous to those 
of Step 2 one may deduce the estimate
 $\log V_k^{(s)}  \le \log V_k^{(0)}+C_k\sqrt{\log V_k^{(s-1)}}$, and hence again 
follows the stabilization of  $V_k^{(s)},  \alpha_{j,k}^{(s)}, E_s$ for all
 $s\ge s^*:=s^*(k)$. \ 
By the very construction the established value $V_k:=V_{\ k}^{(s^*)}$ satisfies 
 {\it all}  inequalities (4.1).

\vspace{1ex}

{\it Step  4.}\ It is easy to check that  all constructions of the Steps 1 -- 3 have been
accomplished by supplementing of {\it minimally required} prime factors to the initial
number $T(p_k)$. Therefore if for some  $\tilde{V}_k$  the assertions of the part 
${\bf (I)}$ of the Theorem 2 hold true, then necessarily  $\tilde{V}_k \ge {V}_k$,
and by virtue of these very inequalities  
(4.1) one has 
$\tilde{\alpha}_{j,k} \ge {\alpha}_{j,k}$ 
for all  $j\le k$,
 i. e.  ${V}_k \ | \ \tilde{V}_k$. 

\vspace{1ex}

{\it Step 5.} The estimates of the part  ${\bf (III)}$ were obtained during the  proof $\ \Box$

\vspace*{\fill}
\clearpage

\ The comparison of Theorems 1 and 2 (cf. (3.3), (4.1)) yields the assertion

\vspace{1ex}

 P r o p o s i t i o n\ \ 2.    
\  {\it In order  the  number   $V_k, \ k\ge 4, $ constructed in Theorem 2,
 to be a one-step  $G$-unimprovable  ($V_k\in {\bf U}_1$) it is necessary 
and sufficient that in addition to (4.1)  {\bf two more} inequalities:} 
 $$\xi(p_k, 0)+\log p_k \ \le \log V_k\ 
\le \xi(p_{k+1}, 0)
\eqno{(4.6)}
$$
{\it hold true.}

\vspace{1ex}

 T a b l e \ 1: { The first 6 one-step 
 $G$-unimprovable numbers computed 
by 

Maple - 13 according to algorithms given in Theorem 2 
with taking into 

account the filtering relationships (4.6). 

The last column contains
the numbers in the estimates of part ${\bf (III)}$.}
 
\vspace{1ex}

\noindent\begin{tabular}{|c|c|l|l|l|}
\hline  $m$ &  $k_m$ & $\qquad \qquad  \quad N_m^*:=V_{k_m}$ 
& $\ G(N^*_m)$ &\ $C_{k_m}$
 \\
\hline  1 & 9 & $ {T}(23)\cdot T(5)\cdot {3}\cdot {2}^3 \ \qquad  \quad \ =\ 160\ 626\ 866\ 400 $  & $1.7374..$ & 1.37..
 \\
\hline 2 & 11&  $ T(31)\cdot T({7})\cdot 3^1
\cdot {2}^4 \qquad \qquad \quad =\ 2.02...\cdot 10^{15} $
&$  1.7368..$ & 1.65..
 \\
\hline 3& 16& $ T(53)\cdot T(7)  \cdot  {3}^2 \cdot {2}^5 \qquad \qquad  \quad =\ 1.97... \cdot 10^{24} $
&$ 1.7434...$ &  1.51..
 \\
\hline 4 & 34 & $ T(139)\cdot T(13)\cdot T(5) \cdot 3^2 \cdot 2^6 \ \quad 
  =\ 5.19... \cdot 10^{63} $
& $1.7582..$ & 1.70..
\\
\hline 5 & 99& $ T(523)\cdot T(29)\cdot T(7)\cdot 5\cdot 3^3\cdot 2^8\  =\ 4.08...\cdot 10^{233} $ &
 $ 1.770728.. $ &  1.67..
\\
\hline 6 & 101 & $ T(547)\cdot T(31)\cdot T(7)\cdot 5 \cdot 3^3\cdot 2^8\  =\ 3.75... \cdot 10^{240} $
&$ 1.770765..$ & 1.78..
 \\
\hline
\end{tabular}

\vspace{2ex}

{\bf 5. \ The infinitude of the set  ${\bf U_1}$.}

\vspace{1ex}

 The key role in proving this assertion play simple  {\it sufficient } 
conditions of {\it $G$-improvability} of the numbers $V_k$
 in terms of the values  $p_k, p_{k+1}, \theta(p_k)$  {\it only}.

\vspace{1ex}

 P r o p o s i t i o n\ \ 3. \ {\bf (i)}  {\it For any   $n>4$ }

 $$ p_{n+1} < \theta(p_n) + 0.5 \sqrt{p_n}\ \Rightarrow 
 \ G(V_np_{n+1})> G(V_n).
\eqno{(5.1)}
$$ 

\vspace{1ex}

{\bf (ii): }  {\it For any $m>4$ }
$$ p_{m} > \theta(p_m) + 4 \sqrt{p_m} \ \Rightarrow  
\ G(V_m/p_{m})> G(V_m).
\eqno{(5.2)}
$$

$\triangle$ \ Indeed, in the case ${\bf (i)}$, taking into account of the
part ${\bf (III)}$  of the Theorem  2 with the {\it lower} $C_k>0.5$  one can conclude that:  
$$ \log V_n>\theta(p_n) + 0.5 \sqrt{p_n} >p_{n+1}>\xi(p_{n+1}, 0),
\eqno{(5.3)} 
$$ 
and then apply  Lemma 1. 

In the case  ${\bf (ii)}$, analogously from the same relationship  ${\bf (III)}$ with  the {\it upper}
estimate  $C_k<4$ we can derive  the chain of inequalities
$$\log V_m< \theta(p_m)+4\sqrt{p_m}
< p_m< \xi(p_m, 0)+\log p_m,
\eqno{(5.4)}
$$
and then take use of Lemma 2 ${\ \Box}$ 
\vspace{1ex}

\vspace*{\fill}
\clearpage

Now we can formulate the main result of this paper.

\vspace{1ex}

T h e o r e m \ \ 3. {\it For any  $M>0$ there exists an integer  $r$ such that $p_r>M$ 
and the number $V_r\ $ is one-step 
 $G$-unimprovable, i. e. } 
$V_r\in {\bf U_1}$.
 
\vspace{1ex}

$\triangle$\  According to the known Littlewood theorem  (1914, cf. [6],  Thm 6.20
and notations on p. (xi)) the difference 
$\psi(x)-x=\Omega_{\pm} (\sqrt{x} \log \log \log x)$.

In particular, from here  it follows  (and that's enough for us), that for arbitrary  $M>100$ 
there are the numbers  $x, y, \ y>2x>M$ such that 
$$x<\theta(x) -10\sqrt{x},  \q y>\theta(y) +10\sqrt{y}
\eqno{(5.5)}
$$
 
Let us denote  $p_n := \max\{p\in \mathbb{P}: p\le x\}, \ p_m := \min \{p\in \mathbb{P}: p\ge y \}$. 

For any  $A\in [-10, 10]$  the function  $t-\theta(t)-A\sqrt{t}$ has the jumps  $\ (-\log p_{k})$ 
 at the points  $p_{k}$ and   {\it increases}  on every semi-segment $(p_k, p_{k+1}], \ p_k>100$; 
therefore  for the numbers  $p_n$ and $p_m$ the inequalities:
$$   p_{n}<\theta(p_n) -8\sqrt{p_n} <
\theta(p_{n-1}) ; \qquad p_{m}>\theta(p_m) +6\sqrt{p_m}.
\eqno{(5.6)}  
$$
hold true. 

Now choose     $r\in \{n-1, \dots , m \}$ such 
that 
$$G(V_r)=\max\{ G(V_k): k\in \{n-1, \dots , m \}  \}.
\eqno{(5.7)}
$$ 

  Proposition  3  jointly  with (5.6) yield the inequalities: 
$$G(V_n)\ge G(V_{n-1}p_n)>G(V_{n-1}); \
 G(V_{m-1})\ge G(V_{m}/p_m)>G(V_{m}).
\eqno{(5.8)}
$$
 
 Hence in turn it follows that   $r$ {\bf is strictly in between } of $n-1$ and  $m$, and thus
 $V_r$ is necessarily  one-step $G$-unimprovable $\ \Box$

\vspace*{\fill}

\clearpage

\vspace{2ex}

\centerline{\bf LIST OF REFERENCES }

\vspace{1ex}

\vspace{1ex}

\noindent
[1]\  Gronwall T.H. {\it Some asymptotic expressions in the theory of numbers}. 

Trans. Amer.  Math. 
Soc. V. 14 (1913), pp. 113-122.

\vspace{1ex}
\noindent
[2]\   Ramanujan S. \ {\it Highly composite numbers, annotated and with a foreword by} 

{\it  J.-L. Nicolas and G. Robin},
\ Ramanujan J. V. 1,  1997, pp. 119 - 153.

\vspace{1ex}

\noindent
[3]\   Robin G. \ {\it Grandes valeurs de la fonction somme des diviseurs et hypoth\`ese de}

{\it  Riemann.} J. Math. Pures Appl. V. 63 (1984), pp. 187 -- 213.

\vspace{1ex}

\noindent
[4]\   Caveney G., Nicolas J.-L., Sondow J. \ {\it Robin's theorem, primes, and a new }

{\it elementary reformulation of Riemann Hypothesis.} INTEGERS 11 (2011),

 \#A33, pp 1-10.

\vspace{1ex}
\noindent
[5]\   Vinogradov\, I. M. {\it Vvedenie v teoriyu chisel} (Russian),
Moscow.:{\ Nauka}, \ 1972.

\vspace{1ex}

\noindent
[6]\  Narkiewicz W. \ {\it The Development of
Prime Number Theory.}

Springer-Verlag Berlin Heidelberg, 2000.

\vspace*{\fill}

\end{document}